\documentclass[a4paper]{article}
\usepackage{latexsym}
\usepackage{amssymb}
\usepackage{caption}

\title{Enumeration of small vertex-transitive \\ directed strongly regular graphs}
\author{A. E. Brouwer \& Š. Gyürki}
\date{\today}

\begin{document}
% 728/92  
\maketitle

\begin{abstract}
We find the 728 vertex-transitive directed strongly regular graphs
with parameters $(v,k,t,\lambda,\mu)$ where $1 \le v \le 39$
and $k < v/2$ and $0 < t < k$.
A table with parameters and automorphism groups and some additional data
is given.
\end{abstract}

\section{Introduction}
A {\em directed strongly regular graph} (DSRG) is a
finite directed graph without loops with adjacency matrix $A$ such that
the linear span of $\{ I, J, A \}$ is closed under matrix multiplication.
Equivalently, a DSRG with parameters $(v,k,t,\lambda,\mu)$
is a directed graph without loops on $v$ vertices, regular with
in-degree and out-degree $k$, such that for any vertex $x$ there are
precisely $t$ vertices $y$ where both $x \to y$ and $y \to x$ are edges,
and for any two distinct vertices $x,y$ where $x \to y$ is an edge
(is not an edge) there are precisely $\lambda$ (resp.~$\mu$)
vertices $z$ such that $x \to z$ and $z \to y$ are both edges.
We have $A^2 = tI + \lambda A + \mu (J-I-A)$.

This concept was introduced by Duval \cite{Duval88}.
A brief discussion, a list with feasible parameter values, known
constructions and a bibliography is given in Brouwer \& Hobart \cite{BH}.

DSRGs come in complementary pairs: together with $A$ also $J-I-A$
satisfies the definition. The complement of a DSRG with parameters
$(v,k,t,\lambda,\mu)$ has parameters $(v,v-k-1,v-2k-1+t,v-2k-2+\mu,v-2k+\lambda)$.
Therefore one can assume $k < v/2$ when classifying DSRGs.
Clearly, $0 \le t \le k$. When $t = k$ the graph is undirected,
a strongly regular graph. When $t = 0$ the graph is a directed tournament
(obtained from a skew Hadamard matrix, cf.~\cite{BH}).
A {\em proper} DSRG is a DSRG with $0 < t < k$.

Fiedler, Klin \& Muzychuk \cite{FKM02} found all DSRGs on at most 20 vertices
with a vertex-transitive automorphism group. Gyürki \cite{Gyurki21} announced
the classification of all vertex-transitive DSRGs on at most 31 vertices.
% Gyürki \cite{Gyurki20} gives pictures
% and a detailed description of the DSRG on at most 12 vertices.
Here we extend this work to cover $v \le 39$.
Both authors did the work independently, with rather different methods
(see Section \ref{methods}), and found the same results.

In this note we enumerate the vertex-transitive DSRGs
with parameters $(v,k,t,\allowbreak\lambda,\mu)$ where $1 \le v \le 39$
and $k < v/2$ and $0 < t < k$ and find that there are 728.
A table with parameters and automorphism groups and some additional data
is given.
The 728 graphs are available at \cite{728graphs}.

\section{Small vertex-transitive DSRGs}

\subsection{Totals}
Below a table giving the total number of graphs up to isomorphism
for each $v \le 39$. The columns headed \#, \#C, \#D give entries
for the number (up to isomorphism) of all vertex-transitive DSRGs
with $k < v/2$, the number of such Cayley graphs,
and the number of such Cayley graphs for the dihedral group.
Given a DSRG $\Gamma$ with adjacency matrix $A$, its {\em reverse}
(obtained by reversing all arrows) has adjacency matrix $A^\top$
and is a DSRG with the same parameters.
The entries in the table are pairs $a/s$, where $a$ (``all'')
is the total count, and $s$ (``symmetric'') is the count of
graphs isomorphic to their reverse.
The number of graphs up to isomorphism and symmetry is $(a + s)/2$.

Note that $v$ cannot be prime (Duval \cite{Duval88}),
and $v$ must be even for \#D.

\captionsetup{justification=centering}
\begin{table}[ht]{\centering
\begin{tabular}{|c|c|c|c|}
\hline
$v$ & \# & \#C & \#D \\
\hline
6 & 1/1 & 1/1 & 1/1 \\
8 & 1/1 & 1/1 & 1/1 \\
10 & 3/1 & 1/1 & 1/1 \\
12 & 6/4 & 6/4 & 3/3 \\
14 & 3/1 & 3/1 & 3/1 \\
16 & 4/2 & 4/2 & 2/2 \\
18 & 17/5 & 12/4 & 7/3 \\
20 & 16/4 & 14/4 & 3/3 \\
21 & 7/1 & 7/1 & -- \\
22 & 7/1 & 7/1 & 7/1 \\
\hline
\end{tabular}~~~~\begin{tabular}{|c|c|c|c|}
\hline
$v$ & \# & \#C & \#D \\
\hline
24 & 58/8 & 50/8 & 7/5 \\
26 & 20/2 & 12/2 & 12/2 \\
27 & 10/2 & 10/2 & -- \\
28 & 12/4 & 9/3 & 9/3 \\
30 & 71/7 & 35/5 & 35/5 \\
32 & 177/25 & 153/25 & 5/3 \\
34 & 47/1 & 33/1 & 33/1 \\
36 & 167/19 & 134/14 & 18/6 \\
38 & 63/1 & 63/1 & 63/1 \\
39 & 38/2 & 38/2 & -- \\
\hline
\end{tabular}
\caption{The total number of vertex-transitive DSRGs with $v \le 39$, $k < v/2$.\\
(There are none of orders 9, 15, 25, 33, 35.)}
}
\end{table}

\subsection{Parameters}
The next table refines the above table, and gives counts
per parameter set.

\begin{center}
{\footnotesize
\begin{tabular}{|c@{~~~}c@{~~~}c@{~~~}c@{~~~}c|c|c|c|}
\hline
$v$ & $k$ & $t$ & $\lambda$ & $\mu$ & \# & \#C & \#D \\
\hline
6 & 2 & 1 & 0 & 1 & 1/1 & 1/1 & 1/1 \\
8 & 3 & 2 & 1 & 1 & 1/1 & 1/1 & 1/1 \\
10 & 4 & 2 & 1 & 2 & 3/1 & 1/1 & 1/1 \\
12 & 3 & 1 & 0 & 1 & 1/1 & 1/1 & 0/0 \\
12 & 4 & 2 & 0 & 2 & 1/1 & 1/1 & 1/1 \\
12 & 5 & 3 & 2 & 2 & 4/2 & 4/2 & 2/2 \\
14 & 6 & 3 & 2 & 3 & 3/1 & 3/1 & 3/1 \\
16 & 7 & 4 & 3 & 3 & 3/1 & 3/1 & 1/1 \\
16 & 7 & 5 & 4 & 2 & 1/1 & 1/1 & 1/1 \\
18 & 4 & 3 & 0 & 1 & 1/1 & 1/1 & 0/0 \\
18 & 5 & 3 & 2 & 1 & 2/0 & 2/0 & 0/0 \\
18 & 6 & 3 & 0 & 3 & 1/1 & 1/1 & 1/1 \\
18 & 7 & 5 & 2 & 3 & 4/0 & 2/0 & 0/0 \\
18 & 8 & 4 & 3 & 4 & 8/2 & 5/1 & 5/1 \\
18 & 8 & 5 & 4 & 3 & 1/1 & 1/1 & 1/1 \\
20 & 4 & 1 & 0 & 1 & 1/1 & 1/1 & 0/0 \\
20 & 7 & 4 & 3 & 2 & 2/0 & 2/0 & 0/0 \\
20 & 8 & 4 & 2 & 4 & 5/1 & 5/1 & 1/1 \\
20 & 9 & 5 & 4 & 4 & 8/2 & 6/2 & 2/2 \\
\hline
\end{tabular}~~\begin{tabular}{|c@{~~~}c@{~~~}c@{~~~}c@{~~~}c|c|c|c|}
\hline
$v$ & $k$ & $t$ & $\lambda$ & $\mu$ & \# & \#C & \#D \\
\hline
21 & 6 & 2 & 1 & 2 & 4/0 & 4/0 & -- \\
21 & 8 & 4 & 3 & 3 & 3/1 & 3/1 & -- \\
22 & 10 & 5 & 4 & 5 & 7/1 & 7/1 & 7/1 \\
24 & 5 & 2 & 1 & 1 & 6/0 & 4/0 & 0/0 \\
24 & 6 & 2 & 0 & 2 & 1/1 & 1/1 & 0/0 \\
24 & 7 & 3 & 2 & 2 & 10/2 & 10/2 & 0/0 \\
24 & 8 & 3 & 2 & 3 & 8/0 & 8/0 & 0/0 \\
24 & 8 & 4 & 0 & 4 & 1/1 & 1/1 & 1/1 \\
24 & 9 & 7 & 2 & 4 & 2/0 & 2/0 & 0/0 \\
24 & 10 & 8 & 4 & 4 & 6/0 & 6/0 & 0/0 \\
24 & 11 & 6 & 5 & 5 & 17/1 & 11/1 & 3/1 \\
24 & 11 & 7 & 6 & 4 & 6/2 & 6/2 & 2/2 \\
24 & 11 & 8 & 7 & 3 & 1/1 & 1/1 & 1/1 \\
26 & 11 & 7 & 4 & 5 & 4/0 & 0/0 & 0/0 \\
26 & 12 & 6 & 5 & 6 & 16/2 & 12/2 & 12/2 \\
27 & 8 & 4 & 3 & 2 & 3/1 & 3/1 & -- \\
27 & 10 & 6 & 3 & 4 & 7/1 & 7/1 & -- \\
28 & 12 & 6 & 4 & 6 & 6/2 & 3/1 & 3/1 \\
28 & 13 & 7 & 6 & 6 & 6/2 & 6/2 & 6/2 \\
\hline
\end{tabular}
}
\end{center}

\begin{center}
\begin{table}[ht]
{\footnotesize
\begin{tabular}{|c@{~~}c@{~~}c@{~~}c@{~~}c|c|c|c|}
\hline
$v$ & $k$ & $t$ & $\lambda$ & $\mu$ & \# & \#C & \#D \\
\hline
30 & 5 & 1 & 0 & 1 & 1/1 & 0/0 & 0/0 \\
30 & 9 & 3 & 2 & 3 & 4/0 & 0/0 & 0/0 \\
30 & 9 & 5 & 4 & 2 & 2/0 & 0/0 & 0/0 \\
30 & 10 & 4 & 2 & 4 & 2/0 & 0/0 & 0/0 \\
30 & 10 & 5 & 0 & 5 & 1/1 & 1/1 & 1/1 \\
30 & 11 & 5 & 4 & 4 & 3/1 & 0/0 & 0/0 \\
30 & 12 & 6 & 3 & 6 & 5/1 & 1/1 & 1/1 \\
30 & 13 & 8 & 5 & 6 & 4/0 & 0/0 & 0/0 \\
30 & 13 & 11 & 6 & 5 & 2/0 & 0/0 & 0/0 \\
30 & 14 & 7 & 6 & 7 & 41/1 & 31/1 & 31/1 \\
30 & 14 & 8 & 7 & 6 & 5/1 & 1/1 & 1/1 \\
30 & 14 & 9 & 8 & 5 & 1/1 & 1/1 & 1/1 \\
32 & 7 & 4 & 3 & 1 & 2/0 & 2/0 & 0/0 \\
32 & 9 & 6 & 1 & 3 & 12/0 & 8/0 & 0/0 \\
32 & 10 & 7 & 3 & 3 & 22/0 & 18/0 & 0/0 \\
32 & 13 & 9 & 4 & 6 & 40/6 & 36/6 & 0/0 \\
32 & 14 & 10 & 6 & 6 & 85/15 & 73/15 & 0/0 \\
32 & 15 & 8 & 7 & 7 & 7/1 & 7/1 & 3/1 \\
32 & 15 & 9 & 8 & 6 & 8/2 & 8/2 & 1/1 \\
32 & 15 & 11 & 10 & 4 & 1/1 & 1/1 & 1/1 \\
\hline
\end{tabular}~~\begin{tabular}{|c@{~~}c@{~~}c@{~~}c@{~~}c|c|c|c|}
\hline
$v$ & $k$ & $t$ & $\lambda$ & $\mu$ & \# & \#C & \#D \\
\hline
34 & 14 & 12 & 5 & 6 & 2/0 & 0/0 & 0/0 \\
34 & 15 & 9 & 6 & 7 & 4/0 & 0/0 & 0/0 \\
34 & 16 & 8 & 7 & 8 & 41/1 & 33/1 & 33/1 \\
36 & 8 & 2 & 1 & 2 & 4/0 & 2/0 & 0/0 \\
36 & 9 & 3 & 0 & 3 & 1/1 & 1/1 & 0/0 \\
36 & 10 & 5 & 2 & 3 & 14/0 & 14/0 & 0/0 \\
36 & 11 & 5 & 4 & 3 & 12/2 & 7/1 & 0/0 \\
36 & 11 & 7 & 6 & 2 & 2/0 & 2/0 & 0/0 \\
36 & 12 & 5 & 2 & 5 & 2/0 & 2/0 & 0/0 \\
36 & 12 & 6 & 0 & 6 & 1/1 & 1/1 & 1/1 \\
36 & 13 & 7 & 4 & 5 & 59/1 & 50/0 & 0/0 \\
36 & 13 & 11 & 2 & 6 & 1/1 & 1/1 & 0/0 \\
36 & 16 & 8 & 6 & 8 & 20/4 & 17/3 & 5/1 \\
36 & 17 & 9 & 8 & 8 & 47/7 & 33/5 & 10/2 \\
36 & 17 & 11 & 10 & 6 & 4/2 & 4/2 & 2/2 \\
38 & 18 & 9 & 8 & 9 & 63/1 & 63/1 & 63/1 \\
39 & 10 & 6 & 1 & 3 & 1/1 & 1/1 & -- \\
39 & 12 & 4 & 3 & 4 & 17/1 & 17/1 & -- \\
39 & 14 & 6 & 5 & 5 & 16/0 & 16/0 & -- \\
39 & 16 & 12 & 7 & 6 & 4/0 & 4/0 & -- \\
\hline
\end{tabular}
}
\caption{Numbers of vertex-transitive DSRGs per parameter set.}
\end{table}
\end{center}

\vspace{-1.5cm}
\subsection{Details}
Finally, we give the automorphism groups and some additional data.
First for the 210 Cayley graphs of dihedral groups, then for the
383 remaining Cayley graphs, then for the remaining 135
vertex-transitive graphs (728 graphs in total). 

\subsubsection{Cayley graphs of dihedral groups}
The columns headed \#, S, Aut, stab, rk give the number of graphs,
whether they are symmetric (isomorphic to their reverse),
(an indication of) the full automorphism group, the point stabilizer,
and the rank of the group action. % , respectively.
For an explanation of the comments, see Section \ref{comments}.

\begin{center}
% \begin{table}[ht]
{\footnotesize
\begin{tabular}{@{~}ccc@{~~~}c@{~~~}c@{~~~}c@{~~~}ccccl}
\hline
\# & S & $v$ & $k$ & $t$ & $\lambda$ & $\mu$ & Aut & stab & rk & comments \\
\hline
% 1x [ 6, 2, 1, 0, 1 ]
1 & S & 6 & 2 & 1 & 0 & 1 & $S_3$ & 1 & 6 & $\pi$-join on $(2,0,0,{-}1,0)/2$ \\

\hline
% 1x [ 8, 3, 2, 1, 1 ]
1 & S & 8 & 3 & 2 & 1 & 1 & $D_8$ & 1 & 8 \\

\hline
% 1x [ 10, 4, 2, 1, 2 ]
1 & S & 10 & 4 & 2 & 1 & 2 & $D_{10}$ & 1 & 10 & Olmez-Song \\ % also $\pi$-join on $(2,0,0,{-}1,0)/2$

\hline
% 1x [ 12, 4, 2, 0, 2 ], 2x [ 12, 5, 3, 2, 2 ]
1 & S & 12 & 4 & 2 & 0 & 2 & $2^6 . S_3$  && 7 & $2$-cocl.~ext of $(6,2,1,0,1)$ \\
1 & S & 12 & 5 & 3 & 2 & 2 & $2^6 . S_3$  && 7 & $2$-cliq.~ext of $(6,2,1,0,1)$ \\
1 & S & 12 & 5 & 3 & 2 & 2 & $D_{12}$ & 1 & 12 \\

\hline
% 3x [ 14, 6, 3, 2, 3 ]
2 & - & 14 & 6 & 3 & 2 & 3 & $7:6$ & 3 & 6 & Olmez-Song \\
1 & S & 14 & 6 & 3 & 2 & 3 & $D_{14}$ & 1 & 14 & Olmez-Song \\ % also $\pi$-join on $(2,0,0,{-}1,0)/2$

\hline
% 1x [ 16, 7, 4, 3, 3 ], 1x [ 16, 7, 5, 4, 2 ]
1 & S & 16 & 7 & 4 & 3 & 3 & $D_{16}$ & 1 & 16 \\
1 & S & 16 & 7 & 5 & 4 & 2 & $2^8 . D_8$  && 9 & $2$-cliq.~ext of $(8,3,2,1,1)$ \\

\hline
% 1x [ 18, 6, 3, 0, 3 ], 5x [ 18, 8, 4, 3, 4 ], 1x [ 18, 8, 5, 4, 3 ]
1 & S & 18 & 6 & 3 & 0 & 3 & ${S_3}^6 . S_3$  && 7 & $3$-cocl.~ext of $(6,2,1,0,1)$ \\
2 & - & 18 & 8 & 4 & 3 & 4 & $(3^3:3):2$ & $3^2$ & 10 & Olmez-Song \\ % also $\pi$-join on $(6,2,1,0,1)/2$ \\
1 & S & 18 & 8 & 4 & 3 & 4 & $D_{18}$ & 1 & 18 & Olmez-Song \\ % also $\pi$-join on $(2,0,0,{-}1,0)/2$ \\
2 & - & 18 & 8 & 4 & 3 & 4 & $D_{18}$ & 1 & 18 & Olmez-Song \\
1 & S & 18 & 8 & 5 & 4 & 3 & ${S_3}^6 . S_3$  && 7 & $3$-cliq.~ext of $(6,2,1,0,1)$ \\

\hline
% 1x [ 20, 8, 4, 2, 4 ], 2x [ 20, 9, 5, 4, 4 ]
1 & S & 20 & 8 & 4 & 2 & 4 & $2^{10} . D_{10}$  && 11 & $2$-cocl.~ext of $(10,4,2,1,2)$ \\
1 & S & 20 & 9 & 5 & 4 & 4 & $2^{10} . D_{10}$  && 11 & $2$-cliq.~ext of $(10,4,2,1,2)$ \\
1 & S & 20 & 9 & 5 & 4 & 4 & $D_{20}$ & 1 & 20 \\

\hline
% 7x [ 22, 10, 5, 4, 5 ]
2 & - & 22 & 10 & 5 & 4 & 5 & $11:10$ & 5 & 6 \\
1 & S & 22 & 10 & 5 & 4 & 5 & $D_{22}$ & 1 & 22 & $\pi$-join on $(2,0,0,{-}1,0)/2$ \\
4 & - & 22 & 10 & 5 & 4 & 5 & $D_{22}$ & 1 & 22 \\

\hline
\multicolumn{11}{r}{\it continued...}
\end{tabular}
}
% \end{table}
\end{center}

% \newpage

\begin{center}
% \begin{table}[ht]
{\footnotesize
\begin{tabular}{@{~}ccc@{~~~}c@{~~~}c@{~~~}c@{~~~}ccccl@{~}}
\hline
\# & S & $v$ & $k$ & $t$ & $\lambda$ & $\mu$ & Aut & stab & rk & comments \\
\hline
% 1x [ 24, 8, 4, 0, 4 ], 3x [ 24, 11, 6, 5, 5 ], 2x [ 24, 11, 7, 6, 4 ], 1x [ 24, 11, 8, 7, 3 ]
1 & S & 24 & 8 & 4 & 0 & 4 & ${S_4}^6 . S_3$  && 7 & $4$-cocl.~ext of $(6,2,1,0,1)$ \\
2 & - & 24 & 11 & 6 & 5 & 5 & $(6 \times S_3):2$ & 3 & 16 \\
1 & S & 24 & 11 & 6 & 5 & 5 & $D_{24}$ & 1 & 24 \\
1 & S & 24 & 11 & 7 & 6 & 4 & ${S_4}^6 . S_3$  && 7 & $4$-cliq.~ext of $(6,2,1,0,1)$ \\
1 & S & 24 & 11 & 7 & 6 & 4 & $2^{12} . D_{12}$  && 13 & $2$-cliq.~ext of $(12,5,3,2,2)$ \\
1 & S & 24 & 11 & 8 & 7 & 3 & ${S_3}^8 . D_8$  && 9 & $3$-cliq.~ext of $(8,3,2,1,1)$ \\

\hline
% 12x [ 26, 12, 6, 5, 6 ]
1 & S & 26 & 12 & 6 & 5 & 6 & $13:6$ & 3 & 10 & eq: 1/1 \\
2 & - & 26 & 12 & 6 & 5 & 6 & $13:6$ & 3 & 10 & Olmez-Song \\
1 & S & 26 & 12 & 6 & 5 & 6 & $D_{26}$ & 1 & 26 & Olmez-Song \\ % also $\pi$-join on $(2,0,0,{-}1,0)/2$
8 & - & 26 & 12 & 6 & 5 & 6 & $D_{26}$ & 1 & 26 & Olmez-Song \\

\hline
% 3x [ 28, 12, 6, 4, 6 ], 6x [ 28, 13, 7, 6, 6 ]
2 & - & 28 & 12 & 6 & 4 & 6 & $2^{14} . (7:6)$  && 7 & $2$-cocl.~ext of $(14,6,3,2,3)$ \\
1 & S & 28 & 12 & 6 & 4 & 6 & $2^{14} . D_{14}$  && 15 & $2$-cocl.~ext of $(14,6,3,2,3)$ \\
2 & - & 28 & 13 & 7 & 6 & 6 & $2^{14} . (7:6)$  && 7 & $2$-cliq.~ext of $(14,6,3,2,3)$ \\
1 & S & 28 & 13 & 7 & 6 & 6 & $2^{14} . D_{14}$  && 15 & $2$-cliq.~ext of $(14,6,3,2,3)$ \\
2 & - & 28 & 13 & 7 & 6 & 6 & $2 \times (7:6)$ & 3 & 12 \\
1 & S & 28 & 13 & 7 & 6 & 6 & $D_{28}$ & 1 & 28 \\

\hline
% 1x [ 30, 10, 5, 0, 5 ], 1x [ 30, 12, 6, 3, 6 ], 31x [ 30, 14, 7, 6, 7 ], 1x [ 30, 14, 8, 7, 6 ], 1x [ 30, 14, 9, 8, 5 ]
1 & S & 30 & 10 & 5 & 0 & 5 & ${S_5}^6 . S_3$  && 7 & $5$-cocl.~ext of $(6,2,1,0,1)$ \\
1 & S & 30 & 12 & 6 & 3 & 6 & ${S_3}^{10} . D_{10}$  && 11 & $3$-cocl.~ext of $(10,4,2,1,2)$ \\
2 & - & 30 & 14 & 7 & 6 & 7 & $3^5:D_{10}$ & $3^4$ & 14 & Olmez-Song \\ % also $\pi$-join on $(6,2,1,0,1)/2$
2 & - & 30 & 14 & 7 & 6 & 7 & $5^3:S_3$ & $5^2$ & 14 & Olmez-Song \\
1 & S & 30 & 14 & 7 & 6 & 7 & $D_{30}$ & 1 & 30 & Olmez-Song \\ % also $\pi$-join on $(2,0,0,{-}1,0)/2$
26 & - & 30 & 14 & 7 & 6 & 7 & $D_{30}$ & 1 & 30 & Olmez-Song \\
1 & S & 30 & 14 & 8 & 7 & 6 & ${S_3}^{10} . D_{10}$ && 11 & $3$-cliq.~ext of $(10,4,2,1,2)$ \\
1 & S & 30 & 14 & 9 & 8 & 5 & ${S_5}^6 . S_3$  && 7 & $5$-cliq.~ext of $(6,2,1,0,1)$ \\

\hline
% 3x [ 32, 15, 8, 7, 7 ], 1x [ 32, 15, 9, 8, 6 ], 1x [ 32, 15, 11, 10, 4 ]
1 & S & 32 & 15 & 8 & 7 & 7 & $D_{32}$ & 1 & 32 \\
2 & - & 32 & 15 & 8 & 7 & 7 & $D_{32}$ & 1 & 32 \\
1 & S & 32 & 15 & 9 & 8 & 6 & $2^{16} . D_{16}$ && 17 & $2$-cliq.~ext of $(16,7,4,3,3)$ \\
1 & S & 32 & 15 & 11 & 10 & 4 & ${S_4}^8 . D_8$ && 9 & $4$-cliq.~ext of $(8,3,2,1,1)$ \\

\hline
% 33x [ 34, 16, 8, 7, 8 ]
1 & S & 34 & 16 & 8 & 7 & 8 & $D_{34}$ & 1 & 34 & Olmez-Song \\ % also $\pi$-join on $(2,0,0,{-}1,0)/2$
30 & - & 34 & 16 & 8 & 7 & 8 & $D_{34}$ & 1 & 34 & Olmez-Song \\
2  & - & 34 & 16 & 8 & 7 & 8 & $D_{34}$ & 1 & 34 & eq: 1/1 \\

\hline
% 1x [ 36, 12, 6, 0, 6 ], 5x [ 36, 16, 8, 6, 8 ], 10x [ 36, 17, 9, 8, 8 ], 2x [ 36, 17, 11, 10, 6 ]
1 & S & 36 & 12 & 6 & 0 & 6 & ${S_6}^6 . S_3$  && 7 & $6$-cocl.~ext of $(6,2,1,0,1)$ \\
2 & - & 36 & 16 & 8 & 6 & 8 & $2^{18} . ((3^3:3):2)$  && 11 & $2$-cocl.~ext of $(18,8,4,3,4)$ \\
1 & S & 36 & 16 & 8 & 6 & 8 & $2^{18} . D_{18}$  && 19 & $2$-cocl.~ext of $(18,8,4,3,4)$ \\
2 & - & 36 & 16 & 8 & 6 & 8 & $2^{18} . D_{18}$  && 19 & $2$-cocl.~ext of $(18,8,4,3,4)$ \\
2 & - & 36 & 17 & 9 & 8 & 8 & $2^{18} . ((3^3:3):2)$  && 11 & $2$-cliq.~ext of $(18,8,4,3,4)$ \\
1 & S & 36 & 17 & 9 & 8 & 8 & $2^{18} . D_{18}$  && 19 & $2$-cliq.~ext of $(18,8,4,3,4)$ \\
2 & - & 36 & 17 & 9 & 8 & 8 & $2^{18} . D_{18}$  && 19 & $2$-cliq.~ext of $(18,8,4,3,4)$ \\
2 & - & 36 & 17 & 9 & 8 & 8 & $2 \times ((3^3:3):2)$ & $3^2$ & 20 \\
1 & S & 36 & 17 & 9 & 8 & 8 & $D_{36}$ & 1 & 36 \\
2 & - & 36 & 17 & 9 & 8 & 8 & $D_{36}$ & 1 & 36 \\
1 & S & 36 & 17 & 11 & 10 & 6 & ${S_6}^6 . S_3$  && 7 & $6$-cliq.~ext of $(6,2,1,0,1)$ \\
1 & S & 36 & 17 & 11 & 10 & 6 & ${S_3}^{12} . D_{12}$  && 13 & $3$-cliq.~ext of $(12,5,3,2,2)$ \\

\hline
% 63x [ 38, 18, 9, 8, 9 ]
2 & - & 38 & 18 & 9 & 8 & 9 & $19:18$ & 9 & 6 & Olmez-Song \\
2 & - & 38 & 18 & 9 & 8 & 9 & $19:6$ & 3 & 14 & Olmez-Song \\
2 & - & 38 & 18 & 9 & 8 & 9 & $19:6$ & 3 & 14 & eq: 1/1 \\
1 & S & 38 & 18 & 9 & 8 & 9 & $D_{38}$ & 1 & 38  & Olmez-Song \\ % also $\pi$-join on $(2,0,0,{-}1,0)/2$
54 & - & 38 & 18 & 9 & 8 & 9 & $D_{38}$ & 1 & 38 & Olmez-Song \\
2 & - & 38 & 18 & 9 & 8 & 9 & $D_{38}$ & 1 & 38 & eq: 1/1 \\
\hline
\end{tabular}
}
% \end{table}
\end{center}

\subsubsection{Further Cayley graphs}

\begin{center}
{\footnotesize
\begin{tabular}{@{~}ccc@{~~~}c@{~~~}c@{~~~}c@{~~~}c@{~~~~~~}c@{~~~~~~}c@{~~~~~~}cl@{~}}
\hline
\# & S & $v$ & $k$ & $t$ & $\lambda$ & $\mu$ & Aut & stab & rk & comments \\
\hline
% 1x [ 12, 3, 1, 0, 1 ]
1 & S & 12 & 3 & 1 & 0 & 1 & $S_4$ & 2 & 7 \\

% 2x [ 12, 5, 3, 2, 2 ]
2 & - & 12 & 5 & 3 & 2 & 2 & $S_4$ & 2 & 7 & Olmez-Song \\

\hline
% 2x [ 16, 7, 4, 3, 3 ]
2 & - & 16 & 7 & 4 & 3 & 3 & $GL(2,3)$ & 3 & 8 \\

\hline
% 1x [ 18, 4, 3, 0, 1 ]
1 & S & 18 & 4 & 3 & 0 & 1 & $(3^2:3):2^2$ & $S_3$ & 7 \\

% 2x [ 18, 5, 3, 2, 1 ]
2 & - & 18 & 5 & 3 & 2 & 1 & $(S_3 \times S_3):2$ & $2^2$ & 8 \\

% 2x [ 18, 7, 5, 2, 3 ]
2 & - & 18 & 7 & 5 & 2 & 3 & $S_3 \times S_3$ & 2 & 10 \\

\hline
\end{tabular}
}
\end{center}

\begin{center}
{\footnotesize
\begin{tabular}{@{~}ccc@{~~~}c@{~~~}c@{~~~}c@{~~~}c@{~~~}c@{~~~}c@{~~}cl@{~}}
\hline
\# & S & $v$ & $k$ & $t$ & $\lambda$ & $\mu$ & Aut & stab & rk & comments \\
\hline
% 1x [ 20, 4, 1, 0, 1 ]
1 & S & 20 & 4 & 1 & 0 & 1 & $S_5$ & $S_3$ & 7 \\

% 2x [ 20, 7, 4, 3, 2 ]
2 & - & 20 & 7 & 4 & 3 & 2 & $S_5$ & $S_3$ & 7 & Olmez-Song \\

% 4x [ 20, 8, 4, 2, 4 ]
2 & - & 20 & 8 & 4 & 2 & 4 & $2^{10} . (5:4)$  && 7 & $2$-cocl.~ext of $(10,4,2,1,2)$ \\
2 & - & 20 & 8 & 4 & 2 & 4 & $5:4$ & 1 & 20 \\

% 4x [ 20, 9, 5, 4, 4 ]
2 & - & 20 & 9 & 5 & 4 & 4 & $2^{10} . (5:4)$  && 7 & $2$-cliq.~ext of $(10,4,2,1,2)$ \\
2 & - & 20 & 9 & 5 & 4 & 4 & $5:4$ & 1 & 20 \\

\hline
% 4x [ 21, 6, 2, 1, 2 ]
2 & - & 21 & 6 & 2 & 1 & 2 & $PSL(3,2)$ & $D_8$ & 6 & Olmez-Song \\
2 & - & 21 & 6 & 2 & 1 & 2 & $7:6$ & 2 & 12 & Olmez-Song \\

% 3x [ 21, 8, 4, 3, 3 ]
1 & S & 21 & 8 & 4 & 3 & 3 & $PSL(3,2)$ & $D_8$ & 6 & Olmez-Song \\
2 & - & 21 & 8 & 4 & 3 & 3 & $7:6$ & 2 & 12 & Olmez-Song \\

\hline
% 4x [ 24, 5, 2, 1, 1 ]
2 & - & 24 & 5 & 2 & 1 & 1 & $2 \times S_4$ & 2 & 16 \\
2 & - & 24 & 5 & 2 & 1 & 1 & $S_4$ & 1 & 24 \\

% 1x [ 24, 6, 2, 0, 2 ]
1 & S & 24 & 6 & 2 & 0 & 2 & $2^{12} . S_4$  && 8 & $2$-cocl.~ext of $(12,3,1,0,1)$ \\

% 10x [ 24, 7, 3, 2, 2 ]
1 & S & 24 & 7 & 3 & 2 & 2 & $2^{12} . S_4$  && 8 & $2$-cliq.~ext of $(12,3,1,0,1)$ \\
1 & S & 24 & 7 & 3 & 2 & 2 & $S_4$ & 1 & 24 \\
8 & - & 24 & 7 & 3 & 2 & 2 & $S_4$ & 1 & 24 \\

% 8x [ 24, 8, 3, 2, 3 ]
8 & - & 24 & 8 & 3 & 2 & 3 & $S_4$ & 1 & 24 & e.g.~$\pi$-join on $(6,2,1,0,1)/3$ \\ % for at least one pair

% 2x [ 24, 9, 7, 2, 4 ]
2 & - & 24 & 9 & 7 & 2 & 4 & $2 \times S_4$ & 2 & 14 \\

% 6x [ 24, 10, 8, 4, 4 ]
6 & - & 24 & 10 & 8 & 4 & 4 & $2 \times S_4$ & 2 & 14 \\

% 8x [ 24, 11, 6, 5, 5 ]
2 & - & 24 & 11 & 6 & 5 & 5 & $2 \times S_4$ & 2 & 16 \\
2 & - & 24 & 11 & 6 & 5 & 5 & $(6 \times 2):2$ & 1 & 24 \\
4 & - & 24 & 11 & 6 & 5 & 5 & $S_4$ & 1 & 24 \\

% 4x [ 24, 11, 7, 6, 4 ]
2 & - & 24 & 11 & 7 & 6 & 4 & $2^{12} . S_4$  && 8 & $2$-cliq.~ext of $(12,5,3,2,2)$ \\
2 & - & 24 & 11 & 7 & 6 & 4 & $2 \times S_4$ & 2 & 14 \\

\hline
% 3x [ 27, 8, 4, 3, 2 ]
1 & S & 27 & 8 & 4 & 3 & 2 & $(3^2 : 3):2^2$ & $2^2$ & 9 \\
2 & - & 27 & 8 & 4 & 3 & 2 & $9:6$ & 2 & 15 \\

% 7x [ 27, 10, 6, 3, 4 ]
1 & S & 27 & 10 & 6 & 3 & 4 & $(3^2 : 3):2^2$ & $2^2$ & 9 \\
2 & - & 27 & 10 & 6 & 3 & 4 & $(3^2 : 3):2^2$ & $2^2$ & 9 \\
4 & - & 27 & 10 & 6 & 3 & 4 & 9:6 & 2 & 15 \\

\hline
% 2x [ 32, 7, 4, 3, 1 ]
2 & - & 32 & 7 & 4 & 3 & 1 & $(S_4 \times S_4):2$ & $S_3 \times S_3$ & 8 \\

% 8x [ 32, 9, 6, 1, 3 ]
2 & - & 32 & 9 & 6 & 1 & 3 & $(((4^2 : 3):2):2$ & $S_3$ & 10 \\
2 & - & 32 & 9 & 6 & 1 & 3 & $(2^3 : 4):2$ & 2 & 20 \\
2 & - & 32 & 9 & 6 & 1 & 3 & $((8:2):2):2$ & 2 & 20 \\
2 & - & 32 & 9 & 6 & 1 & 3 & $2^3 : 4$ & 1 & 32 \\

% 18x [ 32, 10, 7, 3, 3 ]
2 & - & 32 & 10 & 7 & 3 & 3 & $((2^4 : 3):2):2$ & $S_3$ & 10 & \cite{AKMM11}, Prop.~3.5 \\
2 & - & 32 & 10 & 7 & 3 & 3 & $(2^4 : 2):3$ & 3 & 12 \\
2 & - & 32 & 10 & 7 & 3 & 3 & $(4^2 : 2):2$ & 2 & 20 \\
2 & - & 32 & 10 & 7 & 3 & 3 & $((8:2):2):2$ & 2 & 20 & WL-rank is 10 \\
2 & - & 32 & 10 & 7 & 3 & 3 & $2^3 : 4$ & 1 & 32 \\
2 & - & 32 & 10 & 7 & 3 & 3 & $(8:2):2$ & 1 & 32 \\
4 & - & 32 & 10 & 7 & 3 & 3 & $2^4 : 2$ & 1 & 32 \\
2 & - & 32 & 10 & 7 & 3 & 3 & $8 : 2^2$ & 1 & 32 \\

% 36x [ 32, 13, 9, 4, 6 ]
1 & S & 32 & 13 & 9 & 4 & 6 & $(4^2:4):2$ & 4 & 10 \\
2 & - & 32 & 13 & 9 & 4 & 6 & % $(((((2^4 : 3):4):3):2):2):2$
  $[2^9 . 3^2]$  & $S_4 \times S_3$ & 12 \\
2 & - & 32 & 13 & 9 & 4 & 6 & % $(((2 \times ((2 \times (2^4 : 2)):2)):2):2):2$
  $[2^{11}]$ & $2^6$ & 13 \\
2 & - & 32 & 13 & 9 & 4 & 6 & % $(((2 \times ((8 \times 2):2)):2):2):2$
  $[2^9]$ & $2^4$ & 13 \\
2 & - & 32 & 13 & 9 & 4 & 6 & % $(((((8 \times 2):2):2):2):2):2$
  $[2^9]$ & $2 \times D_8$ & 13 \\
2 & - & 32 & 13 & 9 & 4 & 6 & $((2^4 : 2):2):2$ & $2^2$ & 13 \\
1 & S & 32 & 13 & 9 & 4 & 6 & $S_4 \times D_8$ & $S_3$ & 16 \\
2 & - & 32 & 13 & 9 & 4 & 6 & $(2^3 : 4):2$ & 2 & 18 \\
2 & - & 32 & 13 & 9 & 4 & 6 & $(8 : 2^2):2$ & 2 & 18 \\
2 & - & 32 & 13 & 9 & 4 & 6 & $(2^4 : 2):2$ & 2 & 20 \\
2 & - & 32 & 13 & 9 & 4 & 6 & $(2^3 : 4):2$ & 2 & 20 \\
1 & S & 32 & 13 & 9 & 4 & 6 & $((2 \times Q_8):2):2$ & 2 & 24 \\
2 & - & 32 & 13 & 9 & 4 & 6 & $((4 \times 2^2):2):2$ & 2 & 24 \\
2 & - & 32 & 13 & 9 & 4 & 6 & $(8 \times 2):2$ & 1 & 32 \\
4 & - & 32 & 13 & 9 & 4 & 6 & $2^3 : 4$ & 1 & 32 \\
4 & - & 32 & 13 & 9 & 4 & 6 & $(8:2):2$ & 1 & 32 \\
1 & S & 32 & 13 & 9 & 4 & 6 & $4^2 : 2$ & 1 & 32 \\
1 & S & 32 & 13 & 9 & 4 & 6 & $2^4 : 2$ & 1 & 32 \\
1 & S & 32 & 13 & 9 & 4 & 6 & $(4 \times 2^2):2$ & 1 & 32 \\

\hline
\multicolumn{11}{r}{\it continued...}
\end{tabular}
}
% \end{table}
\end{center}

\begin{center}
% \begin{table}[ht]
{\footnotesize
\begin{tabular}{@{~}ccc@{~~~}c@{~~~}c@{~~~}c@{~~~}cc@{~~}c@{~~}cl@{}}
\hline
\# & S & $v$ & $k$ & $t$ & $\lambda$ & $\mu$ & Aut & stab & rk & comments \\
\hline

% 73x [ 32, 14, 10, 6, 6 ]
1 & S & 32 & 14 & 10 & 6 & 6 & $(4^2 : 4):2$ & $2$ & 10 \\
2 & - & 32 & 14 & 10 & 6 & 6 & $[2^{11}]$ & $2^6$ & 13 \\
2 & - & 32 & 14 & 10 & 6 & 6 & $[2^9]$ & $2^4$ & 13 \\
2 & - & 32 & 14 & 10 & 6 & 6 & $[2^9]$ & $2 \times D_8$ & 13 \\
2 & - & 32 & 14 & 10 & 6 & 6 & $[2^7]$ & $2^2$ & 13 \\
2 & - & 32 & 14 & 10 & 6 & 6 & $[2^7]$ & $2^2$ & 14 \\
2 & - & 32 & 14 & 10 & 6 & 6 & $[2^7]$ & $2^2$ & 14 \\
2 & - & 32 & 14 & 10 & 6 & 6 & $[2^9]$ & $2 \times D_8$ & 16 \\
1 & S & 32 & 14 & 10 & 6 & 6 & $(2 \times D_{16}):2$ & 2 & 18 \\
1 & S & 32 & 14 & 10 & 6 & 6 & $(2^2 . (4 \times 2)):2$ & 2 & 18 \\
4 & - & 32 & 14 & 10 & 6 & 6 & $(8 : 2^2):2$ & 2 & 18 \\
8 & - & 32 & 14 & 10 & 6 & 6 & $(2^3 : 4):2$ & 2 & 18 \\
2 & - & 32 & 14 & 10 & 6 & 6 & $[2^7]$ & $2^2$ & 20 \\
2 & - & 32 & 14 & 10 & 6 & 6 & $(2^3 : 4):2$ & 2 & 20 \\
2 & - & 32 & 14 & 10 & 6 & 6 & $(2^4 : 2):2$ & 2 & 20 \\
1 & S & 32 & 14 & 10 & 6 & 6 & $D_8 \times D_8$ & 2 & 24 \\
1 & S & 32 & 14 & 10 & 6 & 6 & $((2 \times Q_8):2):2$ & 2 & 24 \\
2 & - & 32 & 14 & 10 & 6 & 6 & $((4 \times 2^2):2):2$ & 2 & 24 \\
3 & S & 32 & 14 & 10 & 6 & 6 & $4^2 : 2$ & 1 & 32 \\
3 & S & 32 & 14 & 10 & 6 & 6 & $2^4 : 2$ & 1 & 32 \\
2 & S & 32 & 14 & 10 & 6 & 6 & $(4 \times 2^2):2$ & 1 & 32 \\
2 & S & 32 & 14 & 10 & 6 & 6 & $2 \times QD16$ & 1 & 32 \\
2 & - & 32 & 14 & 10 & 6 & 6 & $(8 \times 2):2$ & 1 & 32 \\
8 & - & 32 & 14 & 10 & 6 & 6 & $2^3 : 4$ & 1 & 32 \\
4 & - & 32 & 14 & 10 & 6 & 6 & $(8:2):2$ & 1 & 32 \\
2 & - & 32 & 14 & 10 & 6 & 6 & $2 \times ((4 \times 2):2)$ & 1 & 32 \\
4 & - & 32 & 14 & 10 & 6 & 6 & $2^4 : 2$ & 1 & 32 \\
4 & - & 32 & 14 & 10 & 6 & 6 & $(4 \times 2^2):2$ & 1 & 32 \\

% 4x [ 32, 15, 8, 7, 7 ]
2 & - & 32 & 15 & 8 & 7 & 7 & $QD32$ & 1 & 32 \\
2 & - & 32 & 15 & 8 & 7 & 7 & $SL(2,7):2$ & 7:3 & 8 \\

% 7x [ 32, 15, 9, 8, 6 ]
2 & - & 32 & 15 & 9 & 8 & 6 & $2^{16} . GL(2,3)$  && 9 & $2$-cliq.~ext of $(16,7,4,3,3)$ \\
1 & S & 32 & 15 & 9 & 8 & 6 & $((2^3 : 2^2):3):2$ & $S_3$ & 10 \\
2 & - & 32 & 15 & 9 & 8 & 6 & $(2^3 : 4):2$ & 2 & 20 \\
2 & - & 32 & 15 & 9 & 8 & 6 & $8 : 2^2$ & 1 & 32 \\

\hline
% 2x [ 36, 8, 2, 1, 2 ]
2 & - & 36 & 8 & 2 & 1 & 2 & $3^2 : 4$ & 1 & 36 & eq: 2/4 \\

% 1x [ 36, 9, 3, 0, 3 ]
1 & S & 36 & 9 & 3 & 0 & 3 & ${S_3}^{12} . S_4$  && 8 & $3$-cocl.~ext of $(12,3,1,0,1)$ \\

% 14x [ 36, 10, 5, 2, 3 ]
4 & - & 36 & 10 & 5 & 2 & 3 & $(S_3 \times S_3):2$ & 2 & 21 \\
2 & - & 36 & 10 & 5 & 2 & 3 & $3^2 : 4$ & 1 & 36 \\
8 & - & 36 & 10 & 5 & 2 & 3 & $S_3 \times S_3$ & 1 & 36 \\

% 7x [ 36, 11, 5, 4, 3 ]
1 & S & 36 & 11 & 5 & 4 & 3 & ${S_3}^{12} . S_4$  && 8 & $3$-cliq.~ext of $(12,3,1,0,1)$ \\
2 & - & 36 & 11 & 5 & 4 & 3 & $(3 \times A_4):2$ & 2 & 19 \\
2 & - & 36 & 11 & 5 & 4 & 3 & $3^2 : 4$ & 1 & 36 \\
2 & - & 36 & 11 & 5 & 4 & 3 & $S_3 \times S_3$ & 1 & 36 \\

% 2x [ 36, 11, 7, 6, 2 ]
2 & - & 36 & 11 & 7 & 6 & 2 & $2^{18} . (S_3 \times S_3):2$  && 9 & $2$-cliq.~ext of $(18,5,3,2,1)$ \\

% 2x [ 36, 12, 5, 2, 5 ]
2 & - & 36 & 12 & 5 & 2 & 5 & $(S_3 \times S_3):2$ & 2 & 21 \\

% 50x [ 36, 13, 7, 4, 5 ]
2 & - & 36 & 13 & 7 & 4 & 5 & $[2^7 . 3^5]$ & $[2^5 . 3^3]$ & 9 & $\pi$-join on $(9,4,4,1,2)/3$ \\
4 & - & 36 & 13 & 7 & 4 & 5 & $(3 \times A_4):2$ & 2 & 19 \\
2 & - & 36 & 13 & 7 & 4 & 5 & $A_4 \times S_3$ & 2 & 24 \\
2 & - & 36 & 13 & 7 & 4 & 5 & $3 \times A_4$ & 1 & 36 \\
10 & - & 36 & 13 & 7 & 4 & 5 & $3^2 : 4$ & 1 & 36 \\
30 & - & 36 & 13 & 7 & 4 & 5 & $S_3 \times S_3$ & 1 & 36 \\

% 1x [ 36, 13, 11, 2, 6 ]
1 & S & 36 & 13 & 11 & 2 & 6 & $2 \times ((3 \times A_4):2)$ & $2^2$ & 13 \\

% 12x [ 36, 16, 8, 6, 8 ]
1 & S & 36 & 16 & 8 & 6 & 8 & $2^{18} . (3^2 : 4)$  && 11 & $2$-cocl.~ext of $(18,8,4,3,4)$ \\
2 & - & 36 & 16 & 8 & 6 & 8 & $2^{18} . (3^2:QD16)$  && 7 & $2$-cocl.~ext of $(18,8,4,3,4)$ \\
2 & - & 36 & 16 & 8 & 6 & 8 & $[2^4 . 3^4]$ & $S_3 \times S_3$ & 16 \\
2 & - & 36 & 16 & 8 & 6 & 8 & $(S_3 \times S_3):2$ & 2 & 21 \\
1 & S & 36 & 16 & 8 & 6 & 8 & $(S_3 \times S_3):2$ & 2 & 21 \\
2 & - & 36 & 16 & 8 & 6 & 8 & $3^2 : 4$ & 1 & 36 \\
2 & - & 36 & 16 & 8 & 6 & 8 & $S_3 \times S_3$ & 1 & 36 \\

\hline
\multicolumn{11}{r}{\it continued...}
\end{tabular}
}
% \end{table}
\end{center}

\begin{center}
% \begin{table}[ht]
{\footnotesize
\begin{tabular}{@{~}ccc@{~~~}c@{~~~}c@{~~~}c@{~~~}c@{~~~}c@{~~}c@{~~~}cl@{}}
\hline
\# & S & $v$ & $k$ & $t$ & $\lambda$ & $\mu$ & Aut & stab & rk & comments \\
\hline

% 23x [ 36, 17, 9, 8, 8 ]
1 & S & 36 & 17 & 9 & 8 & 8 & $2^{18} . (3^2 : 4)$  && 11 & $2$-cliq.~ext of $(18,8,4,3,4)$ \\
2 & - & 36 & 17 & 9 & 8 & 8 & $2^{18} . (3^2:QD16)$  && 7 & $2$-cliq.~ext of $(18,8,4,3,4)$ \\
2 & - & 36 & 17 & 9 & 8 & 8 & $[2^4 . 3^4]$ & $S_3 \times S_3$ & 16 \\
4 & - & 36 & 17 & 9 & 8 & 8 & $(3^3 : 3):2^2$ & $3^2$ & 20 \\
1 & S & 36 & 17 & 9 & 8 & 8 & $(S_3 \times S_3):2$ & 2 & 21 \\
2 & - & 36 & 17 & 9 & 8 & 8 & $(S_3 \times S_3):2$ & 2 & 21 \\
2 & - & 36 & 17 & 9 & 8 & 8 & $3^2 : 4$ & 1 & 36 \\
1 & S & 36 & 17 & 9 & 8 & 8 & $S_3 \times S_3$ & 1 & 36 \\
8 & - & 36 & 17 & 9 & 8 & 8 & $S_3 \times S_3$ & 1 & 36 \\

% 2x [ 36, 17, 11, 10, 6 ]
2 & - & 36 & 17 & 11 & 10 & 6 & ${S_3}^{12} . S_4$  && 8 & $3$-cliq.~ext of $(12,5,3,2,2)$ \\

\hline
% 1x [ 39, 10, 6, 1, 3 ]
1 & S & 39 & 10 & 6 & 1 & 3 & $13:6$ & 2 & 21 \\

% 17x [ 39, 12, 4, 3, 4 ]
2 & - & 39 & 12 & 4 & 3 & 4 & $13:12$ & 4 & 12 & Olmez-Song \\
1 & S & 39 & 12 & 4 & 3 & 4 & $13:6$ & 2 & 21 & eq: 1/1 \\
2 & - & 39 & 12 & 4 & 3 & 4 & $13:6$ & 2 & 21 & Olmez-Song \\
12 & - & 39 & 12 & 4 & 3 & 4 & $13:3$ & 1 & 39 & Olmez-Song \\

% 16x [ 39, 14, 6, 5, 5 ]
2 & - & 39 & 14 & 6 & 5 & 5 & $13:12$ & 4 & 12 & Olmez-Song \\
2 & - & 39 & 14 & 6 & 5 & 5 & $13:6$ & 2 & 21 & Olmez-Song \\
12 & - & 39 & 14 & 6 & 5 & 5 & $13:3$ & 1 & 39 & Olmez-Song \\

% 4x [ 39, 16, 12, 7, 6 ]
4 & - & 39 & 16 & 12 & 7 & 6 & $13:6$ & 2 & 21 \\

\hline
\end{tabular}
}
% \end{table}
\end{center}

\subsubsection{Further vertex-transitive graphs}

% Here we list the non-Cayley vertex-transitive graphs.

\begin{center}
% \begin{table}[ht]
{\footnotesize
\begin{tabular}{@{~}ccc@{~~~}c@{~~~}c@{~~~}c@{~~~}cc@{~~}c@{~~}cl@{}}
\hline
\# & S & $v$ & $k$ & $t$ & $\lambda$ & $\mu$ & Aut & stab & rk & comments \\
\hline
% 2x [ 10, 4, 2, 1, 2 ]
2 & - & 10 & 4 & 2 & 1 & 2 & $5:4$ & 2 & 6 & Olmez-Song \\

\hline
% 2x [ 18, 7, 5, 2, 3 ]
2 & - & 18 & 7 & 5 & 2 & 3 & $3^2 : 4$ & 2 & 10 \\

% 3x [ 18, 8, 4, 3, 4 ]
2 & - & 18 & 8 & 4 & 3 & 4 & $3^2:QD16$ & $D_8$ & 6 & Duval \cite{Duval88}, Theorem 5.6 \\
1 & S & 18 & 8 & 4 & 3 & 4 & $3^2 : 4$ & 2 & 10 & eq: 1/1 \\

\hline
% 2x [ 20, 9, 5, 4, 4 ]
2 & - & 20 & 9 & 5 & 4 & 4 & $2 \times (5:4)$ & 2 & 12 \\

\hline
% 2x [ 24, 5, 2, 1, 1 ]
2 & - & 24 & 5 & 2 & 1 & 1 & $A_4:4$ & 2 & 16 \\

% 6x [ 24, 11, 6, 5, 5 ]
2 & - & 24 & 11 & 6 & 5 & 5 & $A_5:4$ & $D_{10}$ & 8 \\
4 & - & 24 & 11 & 6 & 5 & 5 & $(S_3 \times S_3):2$ & 3 & 16 \\

\hline
% 4x [ 26, 11, 7, 4, 5 ]
4 & - & 26 & 11 & 7 & 4 & 5 & $13:4$ & 2 & 14 \\

% 4x [ 26, 12, 6, 5, 6 ]
2 & - & 26 & 12 & 6 & 5 & 6 & $13:12$ & 6 & 6 & Duval \cite{Duval88}, Theorem 5.6 \\
2 & - & 26 & 12 & 6 & 5 & 6 & $13:4$ & 2 & 14 \\

\hline
% 3x [ 28, 12, 6, 4, 6 ]
1 & S & 28 & 12 & 6 & 4 & 6 &  $PSL(3,2)$ & $S_3$ & 7 \\
2 & - & 28 & 12 & 6 & 4 & 6 &  $PSL(3,2)$ & $S_3$ & 7 \\

\hline
% 1x [ 30, 5, 1, 0, 1 ]
1 & S & 30 & 5 & 1 & 0 & 1 & $S_6$ & $S_4$ & 7 \\

% 4x [ 30, 9, 3, 2, 3 ]
2 & - & 30 & 9 & 3 & 2 & 3 & $S_5$ & $2^2$ & 11 \\
2 & - & 30 & 9 & 3 & 2 & 3 & $A_5$ & 2 & 16 \\

% 2x [ 30, 9, 5, 4, 2 ]
2 & - & 30 & 9 & 5 & 4 & 2 & $S_6$ & $S_4$ & 7 & Olmez-Song \\

% 2x [ 30, 10, 4, 2, 4 ]
2 & - & 30 & 10 & 4 & 2 & 4 & $A_5$ & 2 & 16 \\

% 3x [ 30, 11, 5, 4, 4 ]
1 & S & 30 & 11 & 5 & 4 & 4 & $S_5$ & $2^2$ & 11 \\
2 & - & 30 & 11 & 5 & 4 & 4 & $A_5$ & 2 & 16 \\

% 4x [ 30, 12, 6, 3, 6 ]
2 & - & 30 & 12 & 6 & 3 & 6 & ${S_3}^{10} . (5:4)$ && 7 & $3$-cocl.~ext of $(10,4,2,1,2)$ \\
2 & - & 30 & 12 & 6 & 3 & 6 & $S_5$ & $2^2$ & 11 \\

% 4x [ 30, 13, 8, 5, 6 ]
2 & - & 30 & 13 & 8 & 5 & 6 & $5^3:(A_4:4)$ & $2 {\times} D_{10} {\times} D_{10}$ & 10 & $\pi$-join on $(10,3,3,0,1)/2$ \\
2 & - & 30 & 13 & 8 & 5 & 6 & $15:4$ & 2 & 18 \\

% 2x [ 30, 13, 11, 6, 5 ]
2 & - & 30 & 13 & 11 & 6 & 5 & $S_5$ & $2^2$ & 11 \\

% 10x [ 30, 14, 7, 6, 7 ]
2 & - & 30 & 14 & 7 & 6 & 7 & $3^5:(5:4)$ & $3^2 \times (3^2 : 2)$ & 10 & Olmez-Song \\
2 & - & 30 & 14 & 7 & 6 & 7 & $5^3:(A_4:4)$ & $2 {\times} D_{10} {\times} D_{10}$ & 10 & Olmez-Song \\ % also $\pi$-join on $(10,4,2,1,2)/2$ \\
2 & - & 30 & 14 & 7 & 6 & 7 & $5^3:(A_4:4)$ & $2 {\times} D_{10} {\times} D_{10}$ & 10 & $\pi$-join on $(10,4,2,1,2)/2$ \\ % eq: 1/1 \\
4 & - & 30 & 14 & 7 & 6 & 7 & $15:4$ & 2 & 18 & Olmez-Song \\

% 4x [ 30, 14, 8, 7, 6 ]
2 & - & 30 & 14 & 8 & 7 & 6 & ${S_3}^{10} . (5:4)$ && 7 & $3$-cliq.~ext of $(10,4,2,1,2)$ \\
2 & - & 30 & 14 & 8 & 7 & 6 & $A_5$ & 2 & 16 \\

\hline
\multicolumn{11}{r}{\it continued...}
\end{tabular}
}
% \end{table}
\end{center}

\begin{center}
% \begin{table}[ht]
{\footnotesize
\begin{tabular}{@{~}ccc@{~~~}c@{~~~}c@{~~~}c@{~~~}cc@{~~}c@{~~~}cl@{~}}
\hline
\# & S & $v$ & $k$ & $t$ & $\lambda$ & $\mu$ & Aut & stab & rk & comments \\
\hline

% 4x [ 32, 9, 6, 1, 3 ]
4 & - & 32 & 9 & 6 & 1 & 3 & $(2^3 : 4):2$ & 2 & 20 \\

% 4x [ 32, 10, 7, 3, 3 ]
4 & - & 32 & 10 & 7 & 3 & 3 & $(2^3 :4):2$ & 2 & 20 \\

% 4x [ 32, 13, 9, 4, 6 ]
2 & - & 32 & 13 & 9 & 4 & 6 & $[2^9]$ & $2 \times D_8$ & 14 \\
2 & - & 32 & 13 & 9 & 4 & 6 & $(2 \times Q_8):4$ & 2 & 24 \\

% 12x [ 32, 14, 10, 6, 6 ]
2 & - & 32 & 14 & 10 & 6 & 6 & $((16:2):2):2$ & $2^2$ & 14 \\
2 & - & 32 & 14 & 10 & 6 & 6 & $(4.(4 \times 2) : 2):2$ & $2^2$ & 14 \\
2 & - & 32 & 14 & 10 & 6 & 6 & $[2^9]$ & $2 \times D_8$ & 14 \\
4 & - & 32 & 14 & 10 & 6 & 6 & $((2^4 : 2):2):2$ & $2^2$ & 20 \\
2 & - & 32 & 14 & 10 & 6 & 6 & $(2 \times Q_8):4$ & 2 & 24 \\

\hline
% 2x [ 34, 14, 12, 5, 6 ]
2 & - & 34 & 14 & 12 & 5 & 6 & $17:4$ & 2 & 18 \\

% 4x [ 34, 15, 9, 6, 7 ]
4 & - & 34 & 15 & 9 & 6 & 7 & $17:4$ & 2 & 18 \\

% 8x [ 34, 16, 8, 7, 8 ]
2 & - & 34 & 16 & 8 & 7 & 8 & $17:16$ & 8 & 6 & Olmez-Song \\
2 & - & 34 & 16 & 8 & 7 & 8 & $17:8$ & 4 & 10 & Olmez-Song \\
4 & - & 34 & 16 & 8 & 7 & 8 & $17:4$ & 2 & 18 & Olmez-Song \\

\hline
% 2x [ 36, 8, 2, 1, 2 ]
2 & - & 36 & 8 & 2 & 1 & 2 & $((3^2 : Q_8):3):2$ & $D_{12}$ & 7 & Olmez-Song \\

% 5x [ 36, 11, 5, 4, 3 ]
2 & - & 36 & 11 & 5 & 4 & 3 & $((3^2 : Q_8):3):2$ & $D_{12}$ & 7 & Olmez-Song \\
1 & S & 36 & 11 & 5 & 4 & 3 & $3^2 : Q_8$ & 2 & 20 \\
2 & - & 36 & 11 & 5 & 4 & 3 & $3^2 : 8$ & 2 & 20 \\

% 9x [ 36, 13, 7, 4, 5 ]
1 & S & 36 & 13 & 7 & 4 & 5 & $3^2 : Q_8$ & 2 & 20 \\
4 & - & 36 & 13 & 7 & 4 & 5 & $(S_3 \times S_3):2$ & 2 & 20 \\
4 & - & 36 & 13 & 7 & 4 & 5 & $3^2 : 8$ & 2 & 20 \\

% 3x [ 36, 16, 8, 6, 8 ]
2 & - & 36 & 16 & 8 & 6 & 8 & $3^2:QD16$ & 4 & 11 \\
1 & S & 36 & 16 & 8 & 6 & 8 & $(S_3 \times S_3):2$ & 2 & 20 \\

% 14x [ 36, 17, 9, 8, 8 ]
2 & - & 36 & 17 & 9 & 8 & 8 & $3^2:QD16$ & 4 & 11 \\
2 & - & 36 & 17 & 9 & 8 & 8 & $2 \times (3^2:QD16)$ & $D_8$ & 12 \\
2 & - & 36 & 17 & 9 & 8 & 8 & $((3^3:2^2):3):2$ & $3 \times S_3$ & 13 \\
1 & S & 36 & 17 & 9 & 8 & 8 & $2 \times (3^2 : 4)$ & 2 & 20 \\
1 & S & 36 & 17 & 9 & 8 & 8 & $(S_3 \times S_3):2$ & 2 & 20 \\
4 & - & 36 & 17 & 9 & 8 & 8 & $(S_3 \times S_3):2$ & 2 & 20 \\
2 & - & 36 & 17 & 9 & 8 & 8 & $(3^3:2^2):3$ & $3^2$ & 20 \\

\hline
\end{tabular}
}
% \end{table}
\end{center}

\section{Constructions and other comments} \label{comments}

\subsection{(Co)clique extension}
The {\em $m$-coclique extension} of a graph is the graph obtained by
replacing each vertex by a set of $m$ mutually nonadjacent vertices.
Starting with a DSRG with parameters $(v,k,t,\lambda,\mu)$ where $t=\mu$,
one obtains one with parameters $(mv,mk,mt,m\lambda,m\mu)$.

The {\em $m$-clique extension} of a graph is the graph obtained by
replacing each vertex by a set of $m$ mutually adjacent vertices.
Starting with a DSRG with parameters $(v,k,t,\lambda,\mu)$ where $t=\lambda+1$,
one obtains one with parameters $(mv,m(k+1)-1,m(t+1)-1,m(\lambda+2)-2,m\mu)$.

\subsection{Olmez-Song construction}
Let $X$ be a set of size $m = sb+1$. For each $x \in X$ let $\Pi^{(x)}$
be a partition of $X \setminus \{x\}$ into $s$ pieces of size $b$.
We construct DSRGs with parameters
% $(ms,bs,bd,(b-1)d,bd)$ with $d=1$
$(ms,bs,b,b-1,b)$ and $(ms,bs+s-1,b+s-1,b+s-2,b+1)$.
Take as vertex set $V = \{ (x,B) \mid x \in X, B \in \Pi^{(x)}\}$.
Let $(x,B) \to (y,C)$ when $x \in C$ or, in the second case, $x=y$, $B \ne C$.
See Olmez \& Song \cite{OS14}, where more general constructions are given.
If some group $G$ is transitive on $X$ with point stabilizers $G_x$
that are transitive on $\Pi^{(x)}$, then the resulting graphs are
vertex-transitive.
%
% For example, if $s=2$, then then dihedral group generated by $C_m$
% and $z \mapsto -z$ will produce vertex-transitive graphs if
% $\Pi^{(0)} = \{A,-A\}$ (mod $m$).
%
For $s=2$ and multiplier $-1$, the former construction was already given
by Duval \cite{Duval88}, Theorem 6.1.

\subsection{In/out-equivalent vertices}
If $t = \mu$, it can happen that two vertices have the same set of
in-neighbours or the same set of out-neighbours. For example,
in the first Olmez-Song construction vertices come in classes
of size $s$ where all vertices in a class have the same out-neighbours.
A table comment `eq: i/o' gives the class sizes $i, o$ for a graph and its reverse,
mainly to distinguish them from Olmez-Song graphs.

\subsection{$\pi$-join}
Start with a DSRG $\Gamma$ with parameters $(v,k,t,\lambda,\mu)$, where
the vertex set $V$ is partitioned into $a$ parts of size $b$, so that $v = ab$.
Suppose that for any vertex $p$ and any part $P$ the number of arrows
from $P$ to $p$ is $\mu$ if $p \notin P$ and $\lambda+b-k$ if $p \in P$.
Then for each nonnegative integer $u$ there is a DSRG called the $\pi$-join
with parameters $(v,k,t,\lambda,\mu)+u(av,v,b,b,b)$ (Gyürki \cite{Gyurki16}).
It can be constructed as follows:
Take as vertex set $V \times \{0,1,\ldots,ua\}$. Put a copy of $\Gamma$
on each layer $V \times \{i\}$. Between different layers add arrows
$(x,g) \to (y,h)$, where $y$ lies in part~$i$ ($1 \le i \le a$),
when $h \in g+(i-1)u+\{1,\ldots,u\}~({\rm mod}~ua+1)$.

\subsection{Permutation rank and WL-rank}
Given any structure with transitive automorphism group,
the permutation rank is the number of orbits of the point stabilizer.
Given any colored graph, the WL-rank is the number of relations
in the coarsest coherent configuration refining the colored graph.
Where both apply, the WL-rank is not larger than the permutation rank.
In the above tables we gave the permutation ranks. Only in a single
case, namely for one pair of DSRGs with parameters $(32,10,7,3,3)$,
do both ranks differ.

\section{The enumeration process} \label{methods}
Both authors independently did the classification.
The first author used {\sf GAP}~\cite{GAP}
together with the {\sf SmallGrp} and {\sf TransGrp}
packages (\cite{BEO26,Hulpke}).

%All vertex-transitive DSRGs of order at most $20$ were classified
% in~\cite{FKM02}, and this classification was extended up to order $31$
% in \cite{Gyurki21}.
The second author used Hanaki \& Miyamoto's catalogue \cite{HM}
of association schemes combined with the fact that for searching
vertex-transitive DSRGs among unions of relations of an association scheme,
it suffices to consider non-commutative schemes of rank at least~$6$ (see \cite{KMMZ}).
% Since currently this catalogue is complete up to order 39, it was decided to
% find all vertex-transitive DSRGs up to this order. 
All computations were performed in {\sf GAP}~\cite{GAP}, using the packages
{\sf GRAPE}~\cite{GRAPE}, {\sf COCO2P}~\cite{COCO2P} and {\sf SetOrbit}~\cite{SetOrbit}.

% avoid a single word on the next page
\enlargethispage{0.3cm}
For each order~$v$, the enumeration proceeded in two stages. First, we
examined all unions of relations in the vertex-transitive, non-thin,
non-commutative association schemes of order~$v$ and rank at least~$6$.
This produced a set $\mathcal G$ containing all non-Cayley
vertex-transitive DSRGs of order~$v$, as well as every Cayley DSRG whose
automorphism group has order greater than~$v$.
It remained to find the Cayley DSRGs whose full automorphism group has
order exactly~$v$. We therefore tested the possible Cayley connection sets,
using {\sf SetOrbit} to retain only one representative from each orbit
under group automorphisms. This yielded the set $\mathcal H$ of all Cayley
DSRGs of order~$v$. After isomorphism testing, $\mathcal G\cup\mathcal H$
is the complete set of vertex-transitive DSRGs of order~$v$,
$\mathcal G\setminus\mathcal H$ consists of the non-Cayley examples, and
$\mathcal H\setminus\mathcal G$ consists of the Cayley examples having
exactly~$v$ automorphisms.

\section*{Acknowledgements}
The second author acknowledges support from the Slovak Scientific
Grant Agency VEGA under Grants No.\,1/0069/23 and No.\,1/0011/25,
and from the Slovak Research and Development Agency under Contracts
No.~APVV-22-0005 and No.~APVV-23-0076.

\end{document}